\documentclass[12pt,reqno]{amsart}

\usepackage{setspace, hyperref, epigraph}

\theoremstyle{definition}

\DeclareMathOperator{\N}{{\mathbb N}}

\DeclareMathOperator{\Q}{{\mathbb Q}}
\DeclareMathOperator{\R}{{\mathbb R}}

\newcommand\astr{{{}^\ast\!\R}}

\author[T. B.]{Tiziana Bascelli} \address{T. Bascelli, Lyceum
Gymnasium ``F. Corradini'', Thiene, Italy}
\email{tiziana.bascelli@liceocorradini.gov.it}

\author[P. B.]{Piotr B\l{}aszczyk}\address{P. B\l{}aszczyk, Institute
of Mathematics, Pedagogical University of Cracow,
Poland}\email{pb@up.krakow.pl}

\author[A. B.]{Alexandre Borovik} \address{A. Borovik, School of
Mathematics, University of Manchester, Oxford Street, Manchester, M13
9PL, United Kingdom} \email{alexandre.borovik@manchester.ac.uk}

\author[V. K.]{Vladimir Kanovei} \address{V. Kanovei, IPPI, Moscow,
and MIIT, Moscow, Russia}\email{kanovei@googlemail.com}

\author[K. K.]{Karin U. Katz}\address{K. Katz, Department of
Mathematics, Bar Ilan University, Ramat Gan 5290002
Israel}\email{katzmik@math.biu.ac.il}

\author[M. K.]{Mikhail G. Katz}\address{M. Katz, Department of
Mathematics, Bar Ilan University, Ramat Gan 5290002
Israel}\email{katzmik@macs.biu.ac.il}

\author[S. K.]{Semen S. Kutateladze}\address{S. Kutateladze, Sobolev
Institute of Mathematics, Novosibirsk State University, Russia}
\email{sskut@math.nsc.ru}

\author[T. M.]{Thomas McGaffey}\address{T. McGaffey, Rice University,
US}\email{Thomas.McGaffey@sjcd.edu}

\author[D. Sc.]{David M. Schaps}\address{D. Schaps, Department of
Classical Studies, Bar Ilan University, Ramat Gan 5290002 Israel}
\email{dschaps@gmail.com}

\author[D. Sh.]{David Sherry}\address{D. Sherry, Department of
Philosophy, Northern Arizona University, Flagstaff, AZ 86011, US}
\email{David.Sherry@nau.edu}

\begin{document}


\thispagestyle{empty}


\title[Cauchy's infinitesimals and his sum theorem]
{Cauchy's infinitesimals, his sum theorem, and foundational paradigms}

\begin{abstract}
Cauchy's \emph{sum theorem} is a prototype of what is today a basic
result on the convergence of a series of functions in undergraduate
analysis.  We seek to interpret Cauchy's proof, and discuss the
related epistemological questions involved in comparing distinct
interpretive paradigms.  Cauchy's proof is often interpreted in the
modern framework of a Weierstrassian paradigm.  We analyze Cauchy's
proof closely and show that it finds closer proxies in a different
modern framework.

\end{abstract}

\maketitle

\tableofcontents

\section{Introduction}
\label{one}

\epigraph{Infinitesimals were to disappear from mathematical practice
in the face of Weierstrass'~$\varepsilon$ and~$\delta$ notation.
\cite[p.\;208]{Bo86}}

Nearly two centuries after Cauchy first published his controversial
\emph{sum theorem}, historians are still arguing over the nature of
its hypotheses and Cauchy's modification thereof in a subsequent
paper.

\subsection{The sum theorem}
\label{s11}

In 1821 Cauchy presented the sum theorem as follows:

\begin{quote}
When the various terms of series [$u_0+ u_1+ u_2+ \ldots+
u_n+u_{n+1}\ldots$] are functions of the same variable~$x$, continuous
with respect to this variable in the neighborhood of a particular
value for which the series converges, the sum~$s$ of the series is
also a continuous function of~$x$ in the neighborhood of this
particular value.
\end{quote}
This definition is found in his textbook \emph{Cours d'Analyse},
Theorem I in section~6.1; the translation above is from
\cite[p.\;90]{BS}.  Cauchy returned to the sum theorem in an article
\cite{Ca53}, and the matter has been debated by mathematicians and
historians alike ever since.

Cauchy is often claimed to have modified%
\footnote{A debate of long standing concerns the issue of whether
Cauchy \emph{modified} or \emph{clarified} the hypothesis of the sum
theorem of 1853 as compared to 1821.  Our analysis of the 1853 text is
independent of this debate and our main conclusions are compatible
with either view.}
the hypothesis of his sum theorem in 1853 to a stronger hypothesis of
uniform convergence.  What did that modification consist of precisely?
Writes G. Arsac:
\begin{quote}
Assez curieusement, Robinson montre que si l'on interpr\`ete de fa\c
con non standard les hypoth\`eses de l'\'enonc\'e de Cauchy de 1853
que nous venons d'\'etudier, alors cet \'enonc\'e implique la
convergence uniforme.  \cite[p.\,134]{Ar13}
\end{quote}
Our text can be viewed as an extended commentary on Arsac's
\emph{rather curious} observation.

The received scholarship on this issue contains a tension between a
pair of contradictory contentions.  On the one hand, it holds that
Cauchy worked with an Archimedean continuum (A-track for short)
exclusively; on the other, it holds that Cauchy introduced the
condition of uniform convergence in 1853.  These contentions are
contradictory because Cauchy's 1853 condition was not formulated in an
A-track fashion and on the contrary involved infinitesimals.

In more detail, uniform convergence is traditionally formulated over
an Archi\-me\-dean continuum in terms of alternating quantifiers which
are conspicuously absent from Cauchy's text, whereas infinitesimals
are conspicuously present there.  From the A-track viewpoint, the
modication required for uniform convergence focuses on how the
remainder term $r_n(x)$ behaves both (1) as $n$ increases and (2) as a
function of $x$.  This is, however, not what Cauchy does in his 1853
paper.  Namely, Cauchy declares $n$ to be infinite, and requires $r_n$
to be infinitesimal for all $x$ (\emph{toujours}).

An infinitesi\-mal-enriched continuum (B-track for short) enables a
characterisation of uniform convergence by means of enlargement of the
domain, along the lines of what Cauchy suggests by means of his term
\emph{toujours} (\emph{always}) as widely discussed in the literature;
see e.g., \cite{La87}, \cite{Cu88}.  We argue that Cauchy's procedures
as developed in his 1853 text admit better proxies in modern
infinitesimal frameworks%
\footnote{Some historians are fond of recycling the claim that Abraham
Robinson used \emph{model theory} to develop his system with
infinitesimals.  What they tend to overlook is not merely the fact
that an alternative construction of the hyperreals via an ultrapower
requires nothing more than a serious undergraduate course in algebra
(covering the existence of a maximal ideal), but more significantly
the distinction between \emph{procedures} and \emph{foundations} (see
Section~\ref{s12}) which highlights the point that whether one uses
Weierstrass's foundations or Robinson's is of little import,
procedurally speaking.}
than in modern Weierstrassian frameworks limited to an Archimedean
continuum.

Thus, a B-framework allows us to follow Cauchy in presenting a
statement of the sum theorem without quantifier alternations, and also
to be faithful to Cauchy in retaining his infinitesimals rather than
providing \emph{paraphrases} for them in terms of either a variable
subordinate to a universal quantifier, or a sequence tending to zero.

\subsection{Wartofski's challenge}
\label{s12}

Wartofsky launched a challenge to the historians of science.
\cite{Wa76} called for a philosophical analysis of the ontology of a
field like mathematics (as practiced during a particular period) as a
prior condition for a historiography that possesses vision.  We
propose a minor step in this direction by introducing a distinction
between, on the one hand, the \emph{procedures} exploited in the
mathematical practice of that period, and, on the other, the
\emph{ontology} (in the narrow sense) of the mathematical entities
like points, numbers, or functions, as used during that period.  We
use Cauchy's sum theorem as a case study.

We seek an approach to Cauchy that interprets his sum theorem in a way
faithful to his own \emph{procedures}.  We argue that, leaving aside
issues of the \emph{ontology of mathematical entities} and their
justification in terms of one or another \emph{foundational}
framework, a B-track (infinitesimal-enriched) framework provides
better proxies for interpreting Cauchy's \emph{procedures} and
inferential moves than does an A-track (Archi\-medean) framework.
More specifically, we will provide a proxy for Cauchy's proof of the
sum theorem and show that our A-track opponents, who criticize
Cauchy's proof, have to ignore some of Cauchy's procedures.

As I. Grattan-Guinness pointed out, Cauchy's procedures in 1853 are
\emph{difficult to interpret} in an A-framework (see Section~\ref{31}
for more details).  Indeed, nowhere does Cauchy rely on quantifier
alternations; furthermore, he exploits the term \emph{always} to refer
to an extension of the possible inputs to the function, which must now
include infinitesimals in addition to ordinary values.  Meanwhile, an
A-framework can't express uniform convergence without quantifier
alternations, nor countenance any extension beyond the real numbers.
We document strawmanship and crude mathematical errors by A-track
scholars in response to B-track interpretations of Cauchy.

The goal of a Cauchy scholar is not to translate Cauchy's procedures
into contemporary mathematics so as to obtain a correct result, but
rather to understand them on his (Cauchy's) own terms, or as close as
possible to his own terms.  We argue that Cauchy's procedures are best
understood in the sense of infinitesimal mathematics, rather than
paraphrased to fit the \emph{Epsilontik} mode.  Recent articles deals
with Cauchy's foundational stance include \cite{So05}, \cite{Ra12},
\cite{Na14}.

\subsection{A re-evaluation}

This article is part of a re-evaluation of the history of
infinitesimal mathematics.  Such studies have recently appeared in
journals as varied as \emph{Erkenntnis} \cite{13f}, \emph{Foundations
of Science} \cite{17b}, \emph{HOPOS} \cite{16a}, \emph{JGPS}
\cite{17a}, \emph{Notices of the American Mathematical Society}
\cite{13a}, \emph{Perspectives on Science} \cite{13e}, \emph{Studia
Leibnitiana} \cite{14c}.%
\footnote{See \url{http://u.cs.biu.ac.il/~katzmik/infinitesimals.html}
for a more detailed list.}
Abraham Robinson's framework has recently become more visible due to
the activity of some high-profile advocates like Terry Tao; see e.g.,
Tao \cite{Ta14}, Tao--Vu \cite{TV}.  The field has also had its share
of high-profile detractors like Errett Bishop \cite{Bi77} and Alain
Connes \cite{CLS}.  Their critiques were analyzed in \cite{11a},
\cite{13d}, \cite{13c}, and \cite{17i}, \cite{Sa17b}.  Additional
criticisms were voiced by J.\;Earman \cite{Ea75}, K.\;Easwaran
\cite{Ea14}, H.\;M.\;Edwards \cite{Ed07}, G.\;Ferraro \cite{Fe04},
J.\;Grabiner\;\cite{Gr81}, J.\;Gray\;\cite{Gr15}, P.\;Halmos
\cite{Ha85}, H.\;Ishiguro \cite{Is90}, K.\;Schubring \cite{Sch},
Y.\;Sergeyev \cite{Se15}, and D.\;Spalt \cite{Sp02}.  These were dealt
with respectively in \cite{13f}, \cite{14a}, \cite{15b}, \cite{17a},
\cite{12b}, \cite{17d}, \cite{16b}, \cite{16a}, \cite{17e},
\cite{17g}, \cite{11b}.

\section{The revised sum  theorem}
\label{two}

Cauchy denotes the~$n$th partial sum of the series
\[
s(x)=u_0(x)+u_1(x)+u_2(x)+\ldots
\]
by~$s_n(x)$.  He also denotes by~$r_n(x)=s(x)-s_n(x)$ the~$n$th
remainder of the series.  Cauchy's revised sum theorem from
\cite[p.\;456--457]{Ca53}, as reprinted in \cite{Ca00}, states the
following.%
\footnote{Here the equation numbers $(1)$ and $(3)$ are in Cauchy's
text as reprinted in \cite{Ca00}.}
\begin{quote}
\textbf{Th\'eor\`eme 1}.  Si les diff\'erents termes de la s\'erie
\begin{equation*}
\tag{1}
\hfill u_0, u_1, u_2, \ldots, u_n, u_{n+1}, \ldots \hfil
\end{equation*}
sont des fonctions de la variable r\'eelle~$x$, continues, par rapport
\`a cette variable, entre des limites donn\'ees; si, d'ailleurs, la
somme
\begin{equation*}
\tag{3}
u_n+u_{n+1} + \ldots u_{n'-1}
\end{equation*}
devient \emph{toujours} infiniment petite pour des valeurs infiniment
grandes des nombres entiers~$n$ et~$n'>n$, la s\'erie~(1) sera
convergente et la somme~$s$ de la s\'erie sera, entre les limites
donn\'ees, fonction continue de la variable~$x$. \cite[p.\;33]{Ca00}
(emphasis added)
\end{quote}
We will highlight several significant aspects of the procedures
exploited by Cauchy in connection with the sum theorem.  Cauchy
proceeds to present the following argument:
\begin{quote}
Si \`a la s\'erie (1) on substitue la s\'erie (2), l'expression (3),
r\'eduite \`a la somme
\begin{equation*}
\tag{4}
\frac{\sin(n+1)x}{n+1}+\frac{\sin(n+2)x}{n+2}+\ldots+\frac{\sin
n'\,x}{n'}
\end{equation*}
s'\'evanouira pour $x = 0$; mais pour des valeurs de $x$
tr\`es-voisines de z\'ero, par exemple pour $x = \frac{1}{n}$, $n$
\'etant un tr\`es-grand nombre, elle pourra diff\'erer notablement de
z\'ero; et si, en attribuant \`a $n$ une tr\`es-grande valeur, on pose
non-seulement $x = \frac{1}{n}$, mais encore $n' = \infty$, la
somme~(4), ou, ce qui revient au m\^eme, le reste $r_n$ de la
s\'erie~(2) se r\'eduira sensiblement \`a l'int\'egrale
\[
\int_1^\infty\frac{\sin x}{x}dx=
\frac{\pi}{2}-1+\frac{1}{1.2.3}\,\frac{1}{3}-
\frac{1}{1.2.3.4.5}\,\frac{1}{5}+\ldots= 0,6244\ldots
\]
\cite[p.\;457]{Ca53}
\end{quote}
Cauchy's use of the substitution ``$x=\frac{1}{n}$'' is puzzling to a
mathematician trained in the Weierstrassian framework.  We will show
that Robinson's framework provides a more successful interpretation.

\subsection{Infinite numbers}
\label{41}

Cauchy considers two instantiations (which he refers to as
\emph{valeurs}) of integer numbers, denoted~$n$ and~$n'$.  Note that
he refers to them as infinitely large \emph{numbers}, rather than
either \emph{quantities} or \emph{sequences} increasing without bound.


\subsection{Cauchy's~$\varepsilon$}
\label{42}

Cauchy points out that in order to prove the sum theorem, one needs to
show that ``le module de~$r_n$ [soit] inf\'erieur \`a un
nombre~$\varepsilon$ aussi petit que l'on voudra.''  This comment can
be formalized in two different ways, [A] and [B]:
\begin{enumerate}
\item[\hbox{[A]}]
$(\forall\varepsilon\in\mathbb{R})(\exists{}m\in\mathbb{N})(\forall
n\in\mathbb{N}) \left[(\varepsilon>0)\wedge(n>m)
\implies(|r_n^{\phantom{I}}|<\varepsilon)\right]$.
\item[\hbox{[B]}]
$(\forall\varepsilon\in\mathbb{R})\left[(\varepsilon>0)
\rightarrow(|r_n^{\phantom{I}}|<\varepsilon)\right]$.
\end{enumerate}
In neither interpretation does~$\varepsilon$ assume infinitesimal
values.

Here~$n$ is a bound variable in formula [A] in the sense of being
subordinate to the universal quantifier.  Meanwhile,~$n$ is
interpreted as a specific value in [B], i.e., it is a free variable.
For formula [B] to be true,~$n$ must be infinite (and~$r_n$
infinitesimal).  Interpretations [A] and~[B] correspond to the two
tracks compared in Section~\ref{53}.

Clarifying the formulas completely would require adding a quantifier
over the variable~$x$ of~$r_n$, as well as specifying explicitly in
[B] that~$n$ is infinite.

\subsection{When is~$n'=\infty$?}
\label{23b}

Cauchy examines the series%
\footnote{The equation number $(2)$ is in Cauchy's original text.}
\begin{equation*}
\tag{2} 
\sin x+\frac{1}{2}\sin 2x+\frac{1}{3}\sin 3x+\ldots
\end{equation*}
closely related to Abel's purported counterexample to the sum theorem,
namely to the sum theorem as stated in 1821.  The series sums to a
(discontinuous) sawtooth waveform.  Cauchy argues that the series does
not in fact satisfy the hypotheses of his 1853 theorem, as follows.
He considers the difference~$s_{n'}-s_n= u_n+u_{n+1}+\ldots+u_{n'-1}$
which translates into the following sum in the case of the series~(2):%
\footnote{The equation number $(4)$ is in Cauchy's original text.}
\begin{equation*}
\tag{4}
\frac{\sin(n+1)x}{n+1}+\frac{\sin(n+2)x}{n+2}+\ldots+\frac{\sin(n')x}{n'}.
\end{equation*}
Cauchy proceeds to assign the value~$\infty$ to the index~$n'$, by
writing~$n'=\infty$, and points out that the sum~(4) is then precisely
the remainder term~$r_n$.

\subsection{Evaluating at~$\frac{1}{n}$}
\label{44}

Cauchy goes on to evaluate the remainder~$r_n(x)$ at~$x=\frac{1}{n}$
and points out in \cite{Ca53} that the result
\begin{quote}
se r\'eduira sensiblement \`a l'int\'egrale
\[
\int_1^\infty\frac{\sin x}{x}dx=
\frac{\pi}{2}-1+\frac{1}{1.2.3}\,\frac{1}{3}-
\frac{1}{1.2.3.4.5}\,\frac{1}{5}+\ldots= 0,6244\ldots
\]
\end{quote}
Here one can evalutate the integral of~$\frac{\sin x}{x}$ by viewing
the series as an infinite Riemann sum; see \cite[p.\;212]{La89}.
Cauchy concludes that the remainder is not infinitesimal and that
therefore the purported counterexample~(2) in fact does not satisfy
the hypotheses of the sum theorem as formulated in 1853.  For a
discussion of Cauchy's proof of the theorem see Section~\ref{54}.

\section{Methodological issues}
\label{seven}

\subsection{Issues raised by Cauchy's article}
\label{six}

Cauchy's 1853 article raises several thorny issues:

\begin{enumerate}
\item
If one is to interpret Cauchy's new hypothesis as uniform convergence,
how can one explain the absence of any quantifier alternations in the
article (cf.\;Section~\ref{31}, claim~GG4)?
\item
If the term \emph{always} is supposed to signal a transition from the
hypothesis of convergence in 1821 to a hypothesis of uniform
convergence in 1853, what amplification does the term provide, given
that the convergence was already required for every~$x$ in the domain
$[x_0,X]$, and there are no further real inputs at which one could
impose additional conditions?
\item
The traditional reading of convergence, whether ordinary or uniform,
leaves no room for indices given by \emph{infinite} numbers~$n$ and
$n'$ referred to by Cauchy (see Section~\ref{41}).
\item
How can~$|r_n|$ be smaller than every real~$\varepsilon>0$ if the
option of interpreting this inequality in terms of alternating
quantifiers is not available (see Section~\ref{42})?
\item
How can~$r_n(x)$ evaluated at~$x=\frac{1}{n}$ produce the definite
integral~$\int_1^\infty\frac{\sin x}{x}dx$ when the answer should
patently be dependent on~$n$ (see Section~\ref{44})?
\end{enumerate}

Before seeking to answer the questions posed in Section~\ref{six} (see
Section~\ref{eight}), we first address some methodological concerns
and point out possible pitfalls.

\subsection{An ontological disclaimer}
\label{51}

Before answering such questions, one needs to clarify what kind of
answer can be expected when dealing with the work of a mathematician
with no access to modern conceptualisations of mathematical entities
like numbers and their set-theoretic ontology prevalent today.

Such a historiographic concern is ubiquitous.  Thus, what could a
modern historian mean, without falling into a trap of presentism, when
he asserts that Lagrange was the first to define the notion of the
derivative~$f'$ of a real function~$f$?  Such an attribution is
apparently problematic since Lagrange shared neither our set-theoretic
ontology of a punctiform continuum%
\footnote{Historians often use the term \emph{punctiform continuum} to
refer to a continuum made out of points (as for example the
traditional set-theoretic $\R$).  Earlier notions of continuum are
generally thought to be non-punctiform (i.e., not made out of points).
The term \emph{punctiform} also has a technical meaning in topology
unrelated to the above distinction.}
given by the set of real numbers, nor our notion of a function as a
relation of a certain type, again in a set-theoretic context.

Some helpful insights in this area were provided by Benacerraf and
Quine, in terms of a dichotomy of \emph{mathematical procedures and
practice} versus \emph{ontology of mathematical entities} like
numbers.  Thus, Quine wrote:
\begin{quote}
Arithmetic is, in this sense, all there is to number: there is no
saying absolutely what the numbers are; there is only arithmetic.
\cite[p.\;198]{Qu}
\end{quote}
In a related development, \cite{Be65} pointed out that if observer~E
learned that the natural numbers ``are" the Zermelo ordinals~$
\emptyset, \, \{\emptyset\}, \, \{\{\emptyset\}\}, \ldots,~$ while
observer~J learned that they are the von\;Neumann ordinals~$
\emptyset, \, \{\emptyset\}, \, \{ \emptyset, \{\emptyset\}\},
\ldots~$ then, strictly speaking, they are dealing with different
things.  Nevertheless, observer~E's actual mathematical practice (and
the mathematical structures and procedures he may be interested in) is
practically the same as observer~J's.  Hence, different ontologies may
underwrite one and the same practice.

The distinction between practice and ontology can be seen as a partial
response to the challenge to the historians of science launched in
\cite{Wa76}, as discussed in Section~\ref{s12}.

When dealing with Lagrange, Cauchy, or for that matter any
mathematician before 1872, we need to keep in mind that their ontology
of mathematical entities like numbers is distinct from ours.  However,
the \emph{procedures} Lagrange uses can often be interpreted and
clarified in terms of their modern proxies, i.e., procedures used by
modern mathematicians, as in the case of the concept of the
derivative.

The dichotomy of practice, procedures, and inferential moves, on the
one hand, versus foundational issues of ontology of mathematical
entities, on the other, is distinct from the familiar opposition
between intuition and rigor.  In a mathematical framework rigorous up
to modern standards, the syntactic procedures and inferential moves
are just as rigorous as the semantic foundational aspects, whether
set-theoretic or category-theoretic.%
\footnote{Nelson's syntactic recasting of Robinson's framework is a
good illustration, in that the logical procedures in Nelson's
framework are certainly up to modern standards of rigor and are
expressed in the classical Zermelo-Fraenkel set theory (ZFC).  With
respect to an enriched set-theoretic language, infinitesimals in
Nelson's Internal Set Theory (IST) can be found within the real number
system $\R$ itself.  The semantic/ontological issues are handled in an
appendix to \cite{Ne77}, showing Nelson's IST to be a conservative
extension of ZFC.}

\subsection{Felix Klein on two tracks}
\label{62}

In addition to an ontological disclaimer as in Section~\ref{51}, one
needs to bear in mind a point made by Felix Klein in 1908, involving
the non-uniqueness of a conceptual framework for interpreting pre-1872
procedures.

Felix Klein described a rivalry of dual approaches in the following
terms.  Having outlined the developments in real analysis associated
with Weierstrass and his followers, Klein pointed out that
\begin{quote}
The scientific mathematics of today is built upon the series of
developments which we have been outlining.  But an essentially
different conception of infinitesimal calculus has been running
parallel with this [conception] through the centuries
\cite[p.\;214]{Kl08}.
\end{quote}
Such a different conception, according to Klein,
\begin{quote}
harks back to old metaphysical speculations concerning the structure
of the continuum according to which this was made up of [...]
infinitely small parts (ibid.).
\end{quote}
Klein went on to formulate a criterion,%
\footnote{A similar criterion was formulated in
\cite[pp.\;116--117]{Fran}.  For a discussion of the Klein--Fraenkel
criterion see \cite[Section\;6.1]{13c}.}
in terms of the mean value theorem, for what would qualify as a
successful theory of infinitesimals, and concluded:
\begin{quote}
I will not say that progress in this direction is impossible, but it
is true that none of the investigators have achieved anything
positive \cite[p.\;219]{Kl08}.
\end{quote}
Thus, the approach based on notions of infinitesimals is not limited
to ``the work of Fermat, Newton, Leibniz and many others in the 17th
and 18th centuries,'' as implied by Victor J. Katz \cite{Ka14}.
Rather, it was very much a current research topic in Felix Klein's
mind.

Of course, Klein in 1908 had no idea at all of modern infinitesimal
frameworks, such as those of A.\;Robinson (see \cite{Ro61}),
Lawvere--Moerdijk--Reyes--Kock (see \cite{Ko06}), or J. Bell
\cite{Be06}.  What Klein was referring to is the \emph{procedural}
issue of how analysis is to be presented, rather than the
\emph{ontological} issue of a specific realisation of an
infinitesimal-enriched number system in the context of a traditional
set theory.

\subsection{Track A and track B in Leibniz}
\label{53}

The dichotomy mentioned by Klein (see Section~\ref{62}) can be
reformulated in the terminology of dual methodology as follows.  One
finds both A-track (i.e., Archimedean) and B-track (Bernoullian, i.e.,
involving infinitesimals) methodologies in historical authors like
Leibniz, Euler, and Cauchy.  Note that scholars attribute the first
systematic use of infinitesimals as a foundational concept to Johann
Bernoulli.  While Leibniz exploited both infinitesimal methods and
``exhaustion'' methods usually interpreted in the context of an
Archimedean continuum, Bernoulli never wavered from the infinitesimal
methodology.  To note the fact of such systematic use by Bernoulli is
not to say that Bernoulli's foundation is adequate, or that it could
distinguish between manipulations with infinitesimals that produce
only true results and those manipulations that can produce false
results.

In addition to ratios of differentials like~$\frac{dy}{dx}$, Leibniz
also considered ratios of ordinary values which he denoted $(d)y$
and~$(d)x$, so that the ratio~$\frac{(d)y}{(d)x}$ would be what we
call today the derivative.  Here~$dx$ and~$(d)x$ were distinct
entities since Leibniz described them as respectively
\emph{inassignable} and \emph{assignable} in his text \emph{Cum
Prodiisset} \cite{Le01c}:
\begin{quote}
\ldots although we may be content with the assignable quantities
$(d)y$,~$(d)v$,~$(d)z$,~$(d)x$, etc., \ldots{} yet it is plain from
what I have said that, at least in our minds, the unassignables
[\emph{inassignables} in the original Latin]~$dx$ and~$dy$ may be
substituted for them by a method of supposition even in the case when
they are evanescent; \ldots{} (as translated in \cite[p.\;153]{Ch}).
\end{quote}
Echoes of Leibnizian terminology are found in Cauchy when he speaks of
\emph{assignable numbers} while proving various convergence tests, see
e.g., \cite[pp.\;37, 40, 87, 369]{BS}.  An alternative term used by
Cauchy and others is \emph{a given quantity}.  \cite{Fi78} shares the
view that historical evaluations of mathematical analysis in general
and Cauchy's work in particular can profitably be made in the light of
both A-track and B-track approaches.

Leibniz repeatedly asserted that his infinitesimals, when compared
to~$1$, violate what is known today as the Archimedean property, viz.,
Euclid's \emph{Elements}, Definition~V.4; see e.g.,
\cite[p.\;322]{Le95b}, as cited in \cite[p.\;14]{Bos}.

\section{From indivisibles to infinitesimals: reception}

\subsection{Indivisibles banned by the jesuits}
\label{ban}

Indivisibles and infinitesimals were perceived as a theological threat
and opposed on doctrinal grounds in the 17th century.  The opposition
was spearheaded by clerics and more specifically by the jesuits.
Tracts opposing indivisibles were composed by jesuits Paul Guldin,
Mario Bettini, and Andr\'e Tacquet \cite[p.\;291]{Re87}.  P.\;Mancosu
writes:
\begin{quote}
Guldin is taking Cavalieri to be composing the continuum out of
indivisibles, a position rejected by the Aristotelian orthodoxy as
having atomistic implications. \ldots{} Against Cavalieri's
proposition that ``all the lines" and ``all the planes" are magnitudes
- they admit of ratios - Guldin argues that ``all the lines \ldots{}
of both figures are infinite; but an infinite has no proportion or
ratio to another infinite."  \cite[p.\;54]{Ma96}
\end{quote}
Tacquet for his part declared that the method of indivisibles ``makes
war upon geometry to such an extent, that if it is not to destroy it,
it must itself be destroyed.''  \cite[p.\;205]{Fe92},
\cite[p.\;119]{Al14}.

In 1632 (the year Galileo was summoned to stand trial over
heliocentrism) the Society's Revisors General led by Jacob Bidermann
banned teaching indivisibles in their colleges \cite{Fe90},
\cite[p.\;198]{Fe92}.  Indivisibles were placed on the Society's list
of permanently banned doctrines in 1651 \cite{He96}.  
The effects of
anti-indivisible bans were still felt in the 18th century, when most
jesuit mathematicians adhered to the methods of Euclidean geometry, to
the exclusion of the new infinitesimal methods:
\begin{quote}
\ldots le grand nombre des math\'ematiciens de [l'Ordre] resta
jusqu'\`a la fin du XVIII$^e$ si\`ecle profond\'ement attach\'e aux
m\'ethodes euclidiennes.  \cite[p.\;77]{Bo27}
\end{quote}
Echoes of such bans were still heard in the 19th century when jesuit
Moigno, who considered himself a student of Cauchy, wrote:
\begin{quote}
In effect, either these magnitudes, smaller than any \emph{given}
magnitude, still have substance and are divisible, or they are simple
and indivisible: in the first case their existence is a chimera,
since, necessarily greater than their half, their quarter, etc., they
are not actually less than any \emph{given} magnitude; in the second
hypothesis, they are no longer mathematical magnitudes, but take on
this quality, this would renounce the idea of the continuum divisible
to infinity, a necessary and fundamental point of departure of all the
mathematical sciences (ibid., xxiv~f.) (as quoted in
\cite[p.\;456]{Sch}) (emphasis added)
\end{quote}
Moigno saw a contradiction where there is none.  Indeed, an
infinitesimal is smaller than any \emph{assignable}, or \emph{given}
magnitude; and has been so at least since Leibniz, as documented in
Section~\ref{53}.  Thus, infinitesimals need not be less than ``their
half, their quarter, etc.,'' because the latter are not assignable.

In fact, the dichotomy of Moigno's critique is reminiscient of the
dichotomy of a critique of Galileo's indivisibles penned by Moigno's
fellow jesuit Orazio Grassi some three centuries earlier.  P.\;Redondi
summarizes it as follows:
\begin{quote}
As for light - composed according to Galileo of indivisible atoms,
more mathematical than physical - in this case, logical contradictions
arise.  Such indivisible atoms must be finite or infinite.  If they
are finite, mathematical difficulties would arise.  If they are
infinite, one runs into all the paradoxes of the separation to
infinity which had already caused Aristotle to discard the atomist
theory \ldots{} \cite[p.\;196]{Re87}.
\end{quote}
This criticism appeared in the first edition of Grassi's book
\emph{Ratio ponderum librae et simbellae}, published in Paris in 1626.
In fact, this criticism of Grassi's
\begin{quote}
exhumed a discounted argument, copied word-for-word from almost any
scholastic philosophy textbook. \ldots{} The Jesuit mathematician
[Paul] Guldin, great opponent of the geometry of indivisibles, and an
excellent Roman friend of [Orazio] Grassi, must have dissuaded him
from repeating such obvious objections.  Thus the second edition of
the \emph{Ratio}, the Neapolitan edition of 1627, omitted as
superfluous the whole section on indivisibles.  \cite[p.\;197]{Re87}.
\end{quote}
Alas, unlike father Grassi, father Moigno had no Paul Guldin to
dissuade him.

\subsection{Unguru's admonition}
\label{64}

There are several approaches to uniform convergence.  These approaches
are equivalent from the viewpoint of modern mathematics (based on
classical logic):
\begin{enumerate}
\item[\mbox{[A]}] the epsilon-delta approach in the A-track context;
\item[\mbox{[B]}] the approach involving a Bernoullian continuum;
\item[\mbox{[C]}] the sequential approach in the A-track context.
\end{enumerate}

The approach usually found in textbooks is the approach [A] via the
epsilon-delta formulation involving quantifier alternations as
detailed in Section~\ref{34}.  The B-track approach is detailed in
Section~\ref{eight}.  An alternative approach [C] via sequences runs
as follows.  A sequence of functions~$r_n$ as~$n$ runs
through~$\mathbb{N}$ converges uniformly to~$0$ on a domain~$D$ if for
each sequence of inputs~$x_n$ in~$D$, the sequence of
outputs~$r_n(x_n)$ tends to zero.  Giusti pursues approaches [A] and
[C] to Cauchy respectively at \cite[p.\;38]{Gi84} and
\cite[p.\;50]{Gi84} (see Section~\ref{five}).  Meanwhile, S.\;Unguru
admonished:
\begin{quote}
It is \ldots{} a historically unforgiveable sin \ldots{} to assume
wrongly that mathematical equivalence is tantamount to historical
equivalence.  \cite[p.\;783]{Un76}
\end{quote}
When seeking a mathematical interpretation of a historical text we
must consider, in addition to questions of mathematical coherence,
also questions of coherence with the actual procedures adopted by the
historical author.

\subsection{Bottazzini: What kept Cauchy from seeing?}
\label{64b}

The points made in Sections \ref{51} through \ref{53} about ontology
vs procedures, as well as multiple possibilities for modern
interpretive frameworks, seem straightforward enough.  However, the
consequences of paying insufficient attention to such distinctions can
be grave.  Thus, U. Bottazzini opines that
\begin{quote}
the techniques that Abel used were those of Cauchy; his definitions of
continuity and of convergence are the same, as are his use of
infinitesimals in demonstrations.  It is precisely this latter fact
that kept both men from finding the weak point in Cauchy's
demonstration and from seeing that there is another form of
convergence, what is today called the \emph{uniform convergence} of a
series of functions.  \cite[pp.\;115-116]{Bo86} (emphasis in the
original)
\end{quote}
Bottazzini's claim that `it is precisely' `the use of infinitesimals
in demonstrations' that `kept [Cauchy] from finding \ldots{} and from
seeing, etc.' only makes sense in the context of an interpretive
paradigm exclusively limited to track A (see Section~\ref{53}).  Such
an approach views the evolution of 19th century analysis as inevitable
progress toward the foundations of analysis purged of infinitesimals
as established by ``the great triumvirate'' \cite[p.\;298]{Bo49}) of
Cantor, Dedekind, and Weierstrass.  In a similar vein, J.\;Dauben
notes:
\begin{quote}
There is nothing in the language or thought of Leibniz, Euler, or
Cauchy \ldots{} that would make them early
Robinsonians. \cite[p.\;180]{Da88}
\end{quote}
``Early Robinsonians'' perhaps not if this means familiarity with
ultrafilters,%
\footnote{For the role of these in a possible construction of the
hyperreals see Section~\ref{f6}.}
but the extent to which the procedures as found in Leibniz, Euler, and
Cauchy admit straightforward proxies in Robinson's framework does not
emerge clearly from the above comment.

It seems almost as if, in an effort to save Cauchy as a herald for
Weierstrassian arithmetical analysis, some A-track enforcers need
carefully to purge the possibility that Cauchy might also be something
quite different, something generally seen \emph{nowadays} as contrary
to the rational structure of his iconage; e.g., their image seems to
be that this \emph{other} was surely a remnant of the irrational in
his methodology.

\subsection{Fraser \emph{vs} Laugwitz}
\label{s84}

In the abstract of his 1987 article in \emph{Historia Mathematica},
Laugwitz is careful to note that he interprets Cauchy's sum theorem
``with his [i.e., Cauchy's] own concepts":
\begin{quote}
It is shown that the famous so-called errors of Cauchy are correct
theorems when interpreted with his own concepts. \cite[p.\;258]{La87}
\end{quote}
In the same abstract, Laugwitz goes on to emphasize: ``No assumptions
on uniformity or on nonstandard numbers are needed."  Indeed, in
section 7 on pages 264--266, Laugwitz gives a lucid discussion of the
sum theorem in terms of Cauchy's infinitesimals, with not a whiff of
modern number systems.  In particular this section does not mention
the article \cite{SL58}.  In a final section 15 entitled ``Attempts
toward theories of infinitesimals," Laugwitz presents a rather general
discussion, with no specific reference to the sum theorem, of how one
might formalize Cauchyan infinitesimals in modern set-theoretic terms.
A reference to \cite{SL58} appears in this final section only.  Thus,
Laugwitz carefully distinguishes between his analysis of Cauchy's
procedures, on the one hand, and the ontological issues of possible
implementations of infinitesimals in a set-theoretic context, on the
other.

Alas, all of Laugwitz's precautions went for naught.  In 2008, he
became a target of damaging innuendo in the updated version of
\emph{The Dictionary of Scientific Biography}.  Here C.\;Fraser writes
as follows in his article on Cauchy:
\begin{quote}
Laugwitz's thesis is that certain of Cauchy's results that were
criticized by later mathematicians are in fact valid \emph{if one is
willing to accept certain assumptions about Cauchy's understanding and
use of infinitesimals.  These assumptions reflect a theory of analysis
and infinitesimals that was worked out by Laugwitz and Curt Schmieden
during the 1950s}.  \cite[p.\;76]{Fr08} (emphasis added)
\end{quote}
What is particularly striking about Fraser's misrepresentation of
Laugwitz's scholarship is Fraser's failure to recognize the dichotomy
of procedure \emph{vs} ontology.  Fraser repeats the performance in
2015 when he writes:
\begin{quote}
Laugwitz, \ldots{} some two decades following the publication by
Schmieden and him of the $\Omega$-calculus commenced to publish a
series of articles arguing that their non-Archimedean formulation of
analysis is well suited to interpret Cauchy's results on series and
integrals.  \cite[p.\;27]{Fr15}
\end{quote}
What Fraser fails to mention is that Laugwitz specifically separated
his analysis of Cauchy's \emph{procedures} from attempts to account
\emph{ontologically} for Cauchy's infinitesimals in modern terms.

\subsection{The rhetoric of reaction}
\label{65}
G.\;Schubring's disagreement with Laugwitz's interpretation of Cauchy
found expression in the following comment:
\begin{quote}
[Giusti's article] spurred Laugwitz to even more detailed attempts to
banish the error and \emph{confirm that Cauchy had used hyper-real
numbers}.  On this basis, he claims, the errors vanish and the
theorems become correct, or, rather, they always were correct (see
Laugwitz 1990, 21).  \cite[p.\;432]{Sch} (emphasis added)
\end{quote}
Here Schubring is \emph{referring to} the article \cite{La90} though
he is most decidedly not \emph{quoting} it.  In fact, there is no
mention of hyperreals on page 21 in \cite{La90}, contrary to
Schubring's claim.  What one does find here is the following comment:
\begin{quote}
The ``mistakes'' show rather, as \emph{experimenta crucis} that one
must understand Cauchy's terms/definitions [\emph{Begriffe}], in the
spirit of the motto,%
\footnote{Here Laugwitz is referring to Cauchy's motto to the effect
that ``Mon but principal a \'et\'e de concilier la rigueur, dont je
m'\'etais fait une loi dans mon \emph{Cours d'analyse} avec la
simplicit\'e que produit la consideration directe des quantit\'es
infiniment petites.''}
in an infinitesimal-mathematical sense.  \cite[p.\;21]{La90}
(translation ours)
\end{quote}
We fully endorse Laugwitz's comment to the effect that Cauchy's
procedures must be understood in the sense of infinitesimal
mathematics, rather than paraphrased to fit the \emph{Epsilontik}
mode.  Note that we are dealing with an author, namely Laugwitz, who
published Cauchy studies in the leading periodicals \emph{Historia
Mathematica} \cite{La87} and \emph{Archive for History of Exact
Sciences} \cite{La89}.%
\footnote{The fact that Laugwitz had published articles in leading
periodicals does not mean that he couldn't have said something wrong.
However, it does suggest the existence of a strawman aspect of
Schubring's sweeping claims against him.}
The idea that Laugwitz would countenance the claim that Cauchy ``had
used hyper-real numbers'' whereas both the term \emph{hyper-real} and
the relevant construction were not introduced until
\cite[p.\;74]{He48}, strikes us as far-fetched.  Meanwhile, in a
colorful us-against-``them'' passage, J.\;Grabiner opines that
\begin{quote}
[Schubring] effectively rebuts the partisans of nonstandard analysis
who wish to make Cauchy one of \emph{them}, using the work of Cauchy's
disciple the Abb\'e Moigno to argue for Cauchy's own intentions.
\cite[p.\;415]{Gr06}. (emphasis added)
\end{quote}
Schubring indeed uses the work of Moigno in an apparent attempt to
refute infinitesimals.  Jesuit Moigno's confusion on the issue of
infinitesimals was detailed in Section~\ref{ban}.

\newcommand\bash{ 
The Oxford English Dictionary defines the verb
\emph{to demonize} as follows:
\begin{quote}
To portray (a person or thing) as wicked and threatening, (now)
esp.\;in an inaccurate or misrepresentative way. (Now the usual
sense.)
\end{quote}
Schubring's comments on Laugwitz misrepresent Laugwitz's position and
portray him as bent on having his way with regard to Cauchy even at
the expense of attributing to Cauchy things (``hyper-reals'') Cauchy
could not possibly have done.  This is as close as a historian can get
to portraying another historian as \emph{wicked}, in line with the OED
definition.  Our conclusion is that Schubring did not properly analyze
Laugwitz but rather misrepresented him via a strawman reading and
sought to ostracize him through ridicule and demonisation.}
%

Contrary to Schubring's claim, Laugwitz did not attribute 20th century
number systems to Cauchy, but rather argued that modern (B-track in
the terminology of Section~\ref{53}) infinitesimal frameworks provide
better proxies for Cauchy's \emph{procedures} than modern
Weierstrassian (A-track in the terminology of Section~\ref{53})
frameworks.  Laugwitz also sought to understand Cauchy's inferential
moves in terms of their modern proxies.

But Schubring's insights don't even reach a quarter of his student
K.\;Viertel's, who writes:
\begin{quote}
\ldots one can find assertions of Cauchy that contradict a consequent
NSA-interpretation.  With respect to the function $\frac{-m}{x}$
and~$x^{-m}$ and their continuity he remarks: Both become infinite
and, as a consequence, discontinuous when $x=0$. (Bradley, Sandifer
2009, 28). \ldots{} In non-standard analysis, however, one
has~$x^{-1}\simeq\infty$ for all~$x\simeq 0$ which entails the
continuity of the function~$\frac{1}{x}$ for the point zero.  Thus
Cauchy's concept of continuity is not always consistent (or is not
always compatible) with the concept of continuity of nonstandard
analysis.  \cite[Section 4.3.1.1, pp.\;56-57]{Vi14} (translation ours)
\end{quote}
Viertel's apparent intuition that all infinite numbers~$x^{-1}$ are
infinitely close to each other makes one wonder how such intuitions
could be squared with Cauchy's apparently distinct infinite
integers~$n$ and~$n'$ (see Sections~\ref{23b} and \ref{order}).
Viertel's claims concerning the \emph{reciprocal function} in
Robinson's framework are of course incorrect.

Such crude mathematical errors are symptomatic of questionable
priority scales.  Viertel's attempt to ``prove'' that the
\emph{reciprocal function}~$\frac{1}{x}$ is continuous at the origin
in nonstandard analysis is doomed from the start since the
extension~$\R\hookrightarrow\astr$ preserves all the first-order
properties under transfer.%
\footnote{The \emph{transfer principle} is a type of theorem that,
depending on the context, asserts that rules, laws or procedures valid
for a certain number system (or more general mathematical structure),
still apply (i.e., are \emph{transfered}) to an extended number system
(or more general mathematical structure).  Thus, the familiar
extension~$\Q\hookrightarrow\R$ preserves the properties of an ordered
field.  To give a negative example, the
extension~$\R\hookrightarrow\R\cup\{\pm\infty\}$ of the real numbers
to the so-called \emph{extended reals} does not preserve the
properties of an ordered field.  The hyperreal extension
$\R\hookrightarrow\astr$ preserves \emph{all} first-order properties.
For example, the identity~$\sin^2 x+\cos^2 x=1$ remains valid for all
hyperreal~$x$, including infinitesimal and infinite values
of~$x\in\astr$.  In particular, the properties of the reciprocal
function remain the same after it is extended to the hyperreal domain.

Arsac tries to explain nonstandard analysis but he seems to be as
unaware of the transfer principle as Viertel: ``Une fonction continue
au sens habituel est une fonction continue aux points standards, mais
elle ne l'est pas obligatoirement aux points non standard.''
\cite[p.\;133]{Ar13} Contrary to his claim, the natural extension of a
continuous function $f$ will be continuous at \emph{all} hyperreal
points $c$ in the sense of the standard definition
$\forall\epsilon>0\;\exists\delta>0\colon|x-c|<\delta\implies|f(x)-f(c)|<\epsilon$,
by the transfer principle.  Indeed, Arsac confused S-continuity and
continuity\ldots}
Anti-infinitesimal rhetoric by A-track scholars appears to be
\emph{reciprocal} to mathematical competence in matters infinitesimal.

\subsection{Gestalt switch}

The debate over Cauchy's understanding of infinitesimals bears some of
the marks of a paradigm shift in the sense of Thomas Kuhn.  A change
is taking place in the historiography of mathematics, involving the
entry of professional mathematicians into the field.  Mathematicians
bring with them the understanding that the same mathematical object or
structure can be defined and described in a variety of mathematical
languages.  Similarly, the same theorem can be proved (or the same
theory developed) in a multitude of languages and conceptual
frameworks.  

When a mathematician of the past used (and in some cases, invented)
several distinct mathematical languages (as was the case for Cauchy),
it is natural for a modern mathematician to exploit modern
formalisations of these languages in the analysis of his predecessor's
work.  This adds depth and added dimension to our vision of the past,
amounting to a \emph{gestalt switch} as formulated in
\cite[pp.\;150--151]{Ku96}:
\begin{quotation}
[T]he transition between competing paradigms cannot be made a step at
a time, forced by logic and neutral experience. Like a gestalt switch,
it must occur all at once (though not necessarily in an instant) or
not at all.  How, then, are scientists brought to make this
transition? Part of the answer is that they are very often
not. Copernicanism made few converts for almost a century after
Copernicus' death. Newton's work was not generally accepted,
particularly on the Continent, for more than half a century after the
\emph{Principia} appeared.  Priestley never accepted the oxygen
theory, nor Lord Kelvin the electromagnetic theory, and so on. The
difficulties of conversion have often been noted by scientists
themselves.  \ldots{} Max Planck, surveying his own career in his
\emph{Scientific Autobiography}, sadly remarked that ``a new
scientific truth does not triumph by convincing its opponents and
making them see the light, but rather because its opponents eventually
die, and a new generation grows up that is familiar with it.''
\end{quotation}
Here Kuhn is quoting \cite [pp.\;33--34] {Pl50}.  In short, science
(including its historiography) makes progress one funeral at a time.

\section{Interpreting Cauchy's infinitesimal mathematics}
\label{eight}

Cutland et al.  agree with Grattan-Guinness (see Section~\ref{31})
concerning the \emph{difficulty} of interpreting Cauchy's procedures
in an Archi\-medean framework:
\begin{quote}
[Cauchy's] modification of his theorem is anything but clear if we
interpret his conception of the continuum as identical with the
`Weierstrassian' concept.  Abraham Robinson [1966] first discusses
Cauchy's original `theorem' \ldots{} \cite[p.\;376]{Cu88}
\end{quote}

\subsection{Robinson's reading}
\label{s51}

The first B-track interpretation of Cauchy's sum theorem was developed
in \cite[pp.\;271--273]{Ro66}.

Modern versions of the sum theorem require an additional hypothesis to
obtain a mathematically correct result, either (i) in terms of uniform
convergence, or (ii) in terms of equicontinuity of the family of
functions.

Robinson considers both possibilities in the context of Cauchy's 1821
formulation of the sum theorem, which is sufficiently ambiguous to
accomodate various readings.  The last mention of the term
\emph{equicontinuous} occurs in \cite{Ro66} on page 272, line 10 from
bottom.

Meanwhile, starting on page 272, line 2 from bottom, Robinson turns to
an analysis of Cauchy's 1853 formulation of the sum theorem and its
proof.  Robinson writes: ``Over thirty years after the publication of
[the sum theorem], Cauchy returned to the problem [in] (Cauchy
[1853]).'' (ibid.)  Robinson concludes:
\begin{quote}
If we interpret this theorem in the sense of Non-standard Analysis, so
that `infiniment petite' is taken to mean `infinitesimal' and
translate `toujours' by `for all $x$' (and not only `for all standard
$x'$), then the condition introduced by Cauchy \ldots amounts
precisely to uniform convergence in accordance with (i) above.''
\cite[p.\;273]{Ro66}
\end{quote}
Thus, in his analysis of Cauchy's 1853 article, Robinson only
envisions the possibility (i) of uniform convergence.  The reason is
that Cauchy's proof contains hints pointing toward such an
interpretation, while it contains no hints in the direction of (ii)
equicontinuity.

Robinson's interpretation is based on a characterisation of uniform
convergence exploiting a hyperreal extension~$\R\hookrightarrow\astr$;
cf.\;\cite[p.\;62, Theorem~6.3]{Da77}.  Down-to-earth treatments of
the hypereals can be found e.g., in \cite{Li88} and \cite{Pr90}; see
Section~\ref{f6}.

Robinson wrote:
\begin{quote}
\textbf{Theorem 4.6.1} The sequence of standard functions
$\{f_n(x)\}$ converges to the standard function~$f(x)$ uniformly
on~$B\subset T$ [where~$T$ is a metric space] if and only
if~$f(p)\approx f_n(p)$ for all~$p\in{}^\ast\!B$ and for all
infinite~$n$. \cite[p.\;116]{Ro66}
\end{quote}
Here the relation of infinite proximity~$f_n(x)\approx f(x)$ signifies
that the difference~$f_n(x)-f(x)$ is infinitesimal.

\subsection{Constructing a hyperreal field}
\label{f6}

We outline a construction of a hyperreal field~$\astr$.  Let~$\R^{\N}$
denote the ring of sequences of real numbers, with arithmetic
operations defined termwise.  Then we
have
\[
\astr=\R^{\N}\!/\,\text{MAX}
\]
where MAX is the maximal ideal consisting of all ``negligible''
sequences~$(u_n)$, i.e., sequences which vanish for a set of indices
of full measure~$\xi$, namely,~$\xi\big(\{n\in\N\colon
u_n=0\}\big)=1$.

Here~$\xi\colon\mathcal{P}(\N)\to\{0,1\}$ (thus $\xi$ takes only two
values,~$0$ and $1$) is a finitely additive measure taking the
value~$1$ on each cofinite set, where~$\mathcal{P}(\N)$ is the set of
subsets of~$\N$.  

The subset~$\mathcal{F}_\xi\subseteq\mathcal{P}(\N)$ consisting of
sets of full measure is called a \emph{free ultrafilter}.  These
originate with \cite{Ta30}.  The construction of a Bernoullian
continuum outlined here was therefore not available prior to that
date.  Further details on the hyperreals can be found in \cite{17f}.
The educational advantages of the infinitesimal approach are discussed
in~\cite{17h}.

\subsection{An example}

To illustrate the mathematical issue involved, we let~$B=[0,1]$ and
consider the sequence of real functions~$f_n(x)=(1-x)^n$ on the
interval~$B$.  Then the limiting function
\[
f(x) = \begin{cases}
1 &\mbox{if } x = 0 \\ 0 & \mbox{if } x > 0\end{cases}
\]
is discontinuous.  The condition that~$f(x)\approx f_n(x)$ for an
infinite~$n$ is satisfied at all real inputs~$x$, but fails at the
infinitesimal input~$x=\frac{1}{n}$ since
$f_n(\frac{1}{n})=(1-\frac{1}{n})^n\approx \frac{1}{e}$
and~$\frac{1}{e}$ is appreciable:~$\frac{1}{e}\not\approx 0$.  The
limit of the sequence can of course be re-written as the sum of a
series, producing a series of continuous real functions with
discontinuous sum.

A similar phenomenon can be illustrated by means of the family
$f_n(x)=\arctan nx$.

\subsection{Cauchy's B-track procedures}
\label{s54}

We will now interpret Cauchy's \emph{procedures} (as opposed to
\emph{ontology}; see Section~\ref{51}) in such a framework.  The
functions~$f_n$ in Robinson's Theorem~4.6.1 correspond to Cauchy's
partial sums~$s_n$.  Our main thesis, in line with the comment by
Grattan-Guinness cited in Section~\ref{32}, is the following.

\begin{quote}
\textbf{Main Thesis.} Cauchy's term \emph{toujours} suggests
\emph{extending} the possible inputs to the function.
\end{quote}
We interpret this \emph{extension} procedure in terms of the following
hyperreal proxy.  We extend the domain~$B$ of a real function~$f\colon
B\to\mathbb{R}$ to its hyperreal extension~${}^\ast\!B$, obtaining a
hyperreal function~$f\colon{}^\ast\!B\to{}^\ast\mathbb{R}$.  Thus,
if~$B$ is Cauchy's interval~$[x_0,X]$ with, say,~$x_0=0$ (as in the
example of the series~(4) in Section~\ref{23b}), then the hyperreal
extension~${}^\ast\!B={}^*[0,X]$ will contain in particular positive
infinitesimals of the form~$\alpha=\frac{1}{n}$ where~$n$ is an
infinite hyperinteger.

The condition~$f_n(x)\approx f(x)$ is equivalent to the remainder term
$r_n=f(x)-f_n(x)$ being infinitesimal.  This is required to hold
\emph{at all hyperreal inputs}~$x$, serving as a proxy for Cauchy's
1853 modified hypothesis with its implied extension of the domain.

Note that Robinson's Theorem 4.6.1 cited in Section~\ref{s51} contains
no quantifier alternations.  This characterisation of uniform
convergence closely parallels the condition found in Cauchy's text.

The condition~$|r_n|<\varepsilon$ for all real positive~$\varepsilon$
amounts to requiring~$r_n$ (more precisely,~$r_n(x)$) to be
infinitesimal, as per item [B] in Section~\ref{42}.  As Cauchy pointed
out, the condition is not satisfied by the series $\sin
x+\frac{1}{2}\sin 2x+\frac{1}{3}\sin 3x+\ldots$

In the hyperreal framework, the quantity~$r_n(\frac{1}{n})$ for all
infinite~$n$ will possess nonzero standard part (shadow) given
precisely by the integral~$\int_1^\infty\frac{\sin x}{x}dx$, in line
with Cauchy's comment cited in Section~\ref{44} (a more detailed
discussion can be found in \cite{Cl71}).  This answers all the queries
formulated in Section~\ref{six}.  

\subsection{An issue of adequacy and refinement}

Bottazzini opines that
\begin{quote}
the language of infinites and infinitesimals that Cauchy used here
seemed ever more \emph{inadequate} to treat the sophisticated and
complex questions then being posed by analysis\ldots{} The problems
posed by the study of nature, such as those Fourier had faced, now
reappeared everywhere in the most delicate questions of ``pure''
analysis and necessarily led to the elaboration of techniques of
inquiry considerably more \emph{refined} than those that had served
French mathematicians at the beginning of the century.  Infinitesimals
were to disappear from mathematical practice in the face of
Weierstrass'~$\varepsilon$ and~$\delta$ notation, \ldots{}
\cite[p.\;208]{Bo86}.  (emphasis added)
\end{quote}
Are infinitesimal techniques ``inadequate,'' or less ``refined'' than
epsilontic ones, as Bottazzini appears to suggest?  We have argued
that, \emph{procedurally} speaking, the infinitesimal techniques of
the classical masters ranging from Leibniz to Euler to Cauchy are more
robust than is generally recognized by historians whose training is
limited to a Weierstrassian framework.

\section{Grabiner's reading}
\label{three}

In this section we will analyze the comments by J.\;Grabiner related
to Cauchy's alleged confusion in the context of the sum theorem.  Two
subtle mathematical distinctions emerged in the middle of the 19th
century: (i) continuity \emph{vs} uniform continuity, and (ii)
convergence \emph{vs} uniform convergence.  We will show that Grabiner
herself is confused about whether Cauchy's purported error consisted
in his failure to distinguish continuity from uniform continuity, or
his failure to distinguish convergence from uniform convergence.

\subsection{Whose confusion?}
\label{21}

Grabiner alleges that Cauchy's 1821 theorem failed to distinguish
pointwise from uniform convergence, and as a result he left his sum
theorem open to some obvious counterexamples.  Grabiner presents her
case for Cauchy's alleged ``confusion'' in the following terms:
\begin{quote}
\ldots{} in treating series of functions, Cauchy did not distinguish
between pointwise and uniform convergence. The verbal formulations
like ``for all" that are involved in choosing deltas did not
distinguish between ``for any epsilon and for all~$x$" and ``for any
$x$, given any epsilon" [19].  Nor was it at all clear in the 1820's
how much depended on this distinction, since proofs about continuity
and convergence were in themselves so novel. We shall see the same
\emph{confusion} between uniform and point-wise convergence as we turn
now to Cauchy's theory of the derivative.  \cite[p.\;191]{Gr83}
(emphasis added)
\end{quote}
Here the bracketed number [19] refers to her footnote 19 that reads as
follows:
\begin{quote}
19. I. Grattan-Guinness, Development of the Foundations of
Mathematical Analysis from Euler to Riemann, M.\;I.\;T. Press, Cambridge
and London, 1970, p.\;123, puts it well: ``Uniform convergence was
tucked away in the word `always,' with no reference to the variable at
all.''  \cite[p.\;194]{Gr83}
\end{quote}
The views of Grattan-Guinness will be examined in Section~\ref{four}.
We now summarize Grabiner's position as follows:
\begin{enumerate}
\item[G1.]
Cauchy did not distinguish between pointwise convergence and uniform
convergence.
\item[G2.]
Cauchy did not distinguish between ``for any epsilon and for all~$x$"
and ``for any~$x$, given any epsilon.''
\item[G3.]
The failure to distinguish between the quantified clauses in item~(2)
is indicative of a ``confusion between uniform and pointwise
convergence'' in Cauchy's procedures.
\item[G4.]
Cauchy's term \emph{always} is somehow related to uniform convergence.
\item[G5.]
Cauchy's use of the term \emph{always} made no reference to
the variable of the function.
\item[G6.]
No mention of infinitesimals whatsoever is made by Grabiner in
connection with Cauchy's treatment of the sum theorem.
\end{enumerate}

\subsection{Quantifier order}
\label{32}

What does Cauchy's mysterious term \emph{always} signify?  A reader
sufficiently curious about it to trace the sources will discover that
Grattan-Guinness (as cited by Grabiner) is not referring to the 1821
\emph{Cours d'Analyse}, but rather to the article \cite{Ca53}, where
Cauchy
\begin{quote}
stated a new theorem which embodied the \emph{extension} of the
necessary and sufficient conditions for convergence, and a revised
theorem 4.4 to cover uniform convergence \ldots{} \cite[p.\;122]{GG70}.
(emphasis added)
\end{quote}
We will get back to Cauchy's \emph{always} in due time.  We can now
analyze Grabiner's claims G1, G2, and G3 as stated in
Section~\ref{21}.  Her claim~G1 is directly contradicted by
Grattan-Guinness in the passage just cited.  Here it turns out that in
1853 Cauchy \emph{did modify} the original 1821 condition by
\emph{extending} it.

Grabiner's claim~G2 involves terms ``for any'', ``for all'', and
``given any'', each of which expresses the universal quantifier.  Thus
claim~G2 postulates a distinction between a pair of quantifier
clauses, namely,
$(\forall x)(\forall \epsilon)$ and~$(\forall \epsilon)(\forall x)$.
However, these two clauses are logically equivalent.  Indeed, the
order of the quantifiers is significant only when one has alternating
existential and universal quantifiers, contrary to her claim~G2.
Quantifier alternation is irrelevant when only universal quantifiers
are involved.

Grabiner's claim~G3 postulates a ``confusion'' on Cauchy's part but is
based on incorrect premises, and is therefore similarly incorrect.
Grabiner's indictment of Cauchy in \cite[p.\;191]{Gr83} does not stand
up to scrutiny, at least if we are to believe her source, namely
Grattan-Guinness.  Since her claims G5 and G6 are based on
Grattan-Guinness, we will postpone their analysis until
Section~\ref{four}.

\subsection{What does it take to elucidate a distinction?}
\label{22}

Turning now to Grabiner's 1981 book, we find the following:
\begin{quote}
The elucidation of the difference between convergence and uniform
convergence by men like Stokes, Weierstrass, and Cauchy himself was
still more than a decade away. \cite[p.\;12]{Gr81}
\end{quote}
Apparently Cauchy did provide an ``elucidation of the difference
between convergence and uniform convergence,'' after all.  The passage
immediately following is particularly revealing:
\begin{quote}
The verbal formulations of limits and \emph{continuity} used by Cauchy
and Bolzano obscured the distinction between ``for any epsilon, there
is a delta that works for all~$x$" and ``for any epsilon and for
all~$x$, there is a delta."  The only tools for handling such
distinctions were words, and the usual formulation with the word
``always'' suggested ``for all~$x$" as well as ``for any epsilon."
(ibid., pp.\;12--13) (emphasis added)
\end{quote}
This passage immediately follows the passage on p.\;12 on uniform
\emph{convergence}, but it appears to deal instead with uniform
\emph{continuity}.  Indeed, the quantifier clause
$\forall\epsilon\exists\delta\forall x$ that Grabiner refers to
customarily appears in the definition of uniform continuity of a
function rather than uniform convergence of a series of functions.%
\footnote{The symbol~$\delta$ is not used in reference to an integer
tending to infinity.}
Grabiner's position as expressed in the 1981 book involves the
following claims (cf.\;those of her 1983 article enumerated in
Section~\ref{21}):
\begin{enumerate}
\item[G7.]
More than a decade later, Cauchy elucidated the difference between
convergence and uniform convergence;
\item[G8.] 
To handle the distinction between continuity and uniform continuity,
Cauchy needed to handle the quantifier-order distinction between
$\forall\epsilon\exists\delta\forall x$, on the one hand, and
$\forall\epsilon\forall x\exists\delta$, on the other.
\item[G9.]
Cauchy only had the term \emph{always} at his disposal.
\end{enumerate}
We will analyze these in Section~\ref{34}.

\subsection{Cauchy's definition of continuity}
\label{toujours}

Having mentioned both continuity and convergence in Section~\ref{22},
we should point out that, while Cauchy's definition of continuity in
1821 was ambiguous to the extent that it is not clear whether he was
defining ordinary or uniform continuity, what is clear is that he
exploited the term \emph{toujours} to refer to the possible inputs to
the function.  He formulated the definition as follows in 1853:
\begin{quote}
\ldots{} une fonction~$u$ de la variable r\'eelle~$x$ sera
\emph{continue}, entre deux limites donn\'ees de~$x$, si, cette
fonction admettant pour chaque valeur interm\'ediaire de~$x$ une
valeur unique et finie, un accroissement infiniment petit attribu\'e
\`a la variable produit \emph{toujours}, entre les limites dont il
s'agit, un accroissement infiniment petit de la fonction
elle-m\^eme. \cite[pp.\;455--456]{Ca53} (emphasis on \emph{continue}
in the original; emphasis on \emph{toujours} ours)
\end{quote}
The phrase ``toujours, entre les limites dont il s'agit'' suggests
that the term \emph{toujours} refers to the possible inputs to the
function.  While Cauchy apparently did not distinguish between
distinct notions of continuity (such as modern notions of pointwise
continuity as opposed to uniform continuity), he did use the term
\emph{toujours} in reference to the possible inputs.

One can surmise that when he uses the term \emph{toujours} in the
definition of convergence in the same article, he is also referring to
the possible inputs.

Cauchy's definition of convergence in 1821 cited in Section\;\ref{s11}
did not use the term \emph{toujours}.  Thus Cauchy was possibly
referring to ordinary convergence.  By contrast, in 1853 the
definition did use the term \emph{toujours} in defining convergence,
which arguably corresponds to uniform convergence (see
Section~\ref{eight} for a more detailed discussion).

\subsection{Uniform convergence today}
\label{34}

To fix notation, recall that today a series~$\sum_{n=0}^\infty u_n(x)$
is said to converge uniformly to a function~$s(x)$ on a domain~$I$ if
\[
(\forall\epsilon>0) (\exists N\in\mathbb{N}) (\forall x\in I) (\forall
m\in\mathbb{N})\left[(m>N)\Rightarrow \left|s(x)-\sum_{n=0}^m
u_n\right|<\epsilon\right]
\]
where all variables and functions are real.  Grabiner's claim~G8 to
the effect that Cauchy failed to distinguish
$\forall\epsilon\exists\delta\forall x$ from~$\forall\epsilon\forall
x\exists\delta$ (see Section~\ref{22} above) is puzzling coming as it
does in the midst of a discussion of the dichotomy of convergence
\emph{vs} uniform convergence of series of functions, rather than
continuity \emph{vs} uniform continuity.  Furthermore, we are not
familiar with any text where Cauchy dealt explicitly with a
distinction related to continuity \emph{vs} uniform continuity.  In
fact, Grabiner herself writes:
\begin{quote}
Weierstrass and his school distinguished--as Cauchy had not--between
pointwise and uniform continuity. \; \; \; \cite[p.\;97]{Gr81}
\end{quote}
We will therefore assume that what Grabiner meant on page 13 of her
book was actually the following modification of her claim~G8:
\begin{enumerate}
\item[G8$'$.] To handle the distinction between \emph{convergence} and
\emph{uniform convergence}, Cauchy needed to handle the
quantifier-order distinction between~$\forall\epsilon\exists N \forall
x$, on the one hand, and~$\forall\epsilon\forall x\exists N$, on the
other.
\end{enumerate}
By page 140 Grabiner seems to have changed her mind once again about
the status of convergence in Cauchy, as she writes:
\begin{quote}
Actually, [Cauchy's] proof implicitly assumed the function to be
uniformly continuous, though he did not distinguish between continuity
and uniform continuity, \emph{just as he had not distinguished between
convergence and uniform convergence}. \cite[p.\;140]{Gr81} (emphasis
added)
\end{quote}
We will assume that she is referring to \cite{Ca21} here rather than
\cite{Ca53}.  Grabiner does not explain Cauchy's use of the term
\emph{always}, but she assumes in claims G7--G9 that the use of this
term somehow corresponds to requiring the stronger \emph{uniform
convergence}.  She also assumes in G8$'$ that successfully dealing
with the dichotomy of convergence \emph{vs} uniform convergence
necessarily requires a quantifier-order distinction, indicative of an
A-track \emph{parti pris}.

\subsection{A-tracking Cauchy's infinitesimals}

Grabiner wrote in 1983:
\begin{quote}
Now we come to the last stage in our chronological list: definition.
In 1823, Cauchy defined the derivative of~$f(x)$ as the limit, when it
exists, of the quotient of differences~$\frac{f(x+h)-f(x)}{h}$ as~$h$
goes to zero [4, pp.\;22--23].  \cite[p.\;204]{Gr83b}
\end{quote}
Her reference ``4, pp.\;22--23'' is to \cite{Ca23}.  Here is Cauchy's
definition that Grabiner claims to report on:
\begin{quote}
Lorsque la fonction~$y=f(x)$ reste continue entre deux limites
donn\'ees de la variable~$x$, et que l'on assigne \`a cette variable
une valeur comprise entre les deux limites dont il s'agit, un
accroissement \emph{infiniment petit}, attribu\'e \`a la variable,
produit un accroissement \emph{infiniment petit} de la fonction
elle-m\^eme.  Par cons\'equent, si l'on pose alors~$\Delta x=i$, les
deux termes du rapport aux diff\'erences
\[
\frac{\Delta y}{\Delta x} = \frac{f(x+i)-f(x)}{i}
\]
seront des quantit\'es \emph{infiniment petites}.  
\cite{Ca23} (emphasis added)
\end{quote}
Notice that the infinitely small have been mentioned \emph{three
times}.  Cauchy continues:
\begin{quote}
Mais, tandis que ces deux termes s'approcheront ind\'efini\-ment et
simultan\'ement de la limite z\'ero, le rapport lui-m\^eme pourra
converger vers une autre limite, soit positive, soit n\'egative.
Cette limite, lorsqu'elle existe, a une valeur d\'etermin\'ee, pour
chaque valeur particuli\`ere de~$x$ \ldots{} on donne \`a la nouvelle
fonction le nom de fonction d\'eriv\'ee, \ldots{} (ibid)
\end{quote}
Grabiner's paraphrase of Cauchy's definition is a blatant example of
failing to see what is right before one's eyes; for Cauchy's
\emph{infinitely small}~$i$ have been systematically purged by
Grabiner, only to be replaced by a post-Weierstrassian imposter
``as~$h$ goes to zero.''  The difference between the Cauchyan and
Weierstrassian notions of limit wsa recently analyzed in \cite{Na14}.

\section{Interpretation by Grattan-Guinness}
\label{four}

As we saw in Section~\ref{three}, Grabiner restricts her discussion of
the sum theorem to mentioning that in 1853, Cauchy introduced uniform
convergence (which is a sufficient condition for the validity of the
sum theorem as stated today).  Grattan-Guinness treats the sum theorem
in more detail.  He notes that
\begin{quote}
Cauchy lived through the appearance of modes of convergence and
returned to his [sum theorem] in a short paper of 1853\ldots{}
[Cauchy] admitted that
\[
\frac{\pi-x}{2}=\sin x+\frac{1}{2}\sin 2x+\frac{1}{3}\sin 3x+\cdots
\]
\ldots{} was\, a counterexample\, with its\, discontinuity of
magnitude $\pi$ when~$x$ equals a multiple of~$\pi$. \; \; \; \; \; \;
\; \; \cite[p.\;122]{GG70}.
\end{quote}

\subsection{Difficult to interpret}
\label{31}

Recall that Cauchy denotes by~$s_n(x)$ the~$n$th partial sum of the
series $s(x)=u_0(x)+u_1(x)+\ldots{}$.  He also denotes
by~$r_n(x)=s(x)-s_n(x)$ the~$n$th remainder of the series.
Grattan-Guinness writes:
\begin{quote}
So a repair was needed, precisely at the point in his original proof
of [the sum theorem] where he had said that
$[r_n(x_0+\alpha)-r_n(x_0)]$ ``becomes insensible at the same time''
as~$r_n(x_0)$.  (ibid.)
\end{quote}
Here an understanding of the term~$r_n(x_0+\alpha)$ is crucial to
interpreting Cauchy's proof, and we will return to it in
Section~\ref{s42}.  Grattan-Guinness proceeds to make the following key
comment:
\begin{quote}
This remark [of Cauchy's] is \emph{difficult to interpret} against
[i.e., in the context of] the classification of modes of uniform
convergence given here \ldots{} since~$\alpha$ is an infinitesimally
small increment of~$x$.  (ibid.) (emphasis added)
\end{quote}
The ``modes of convergence'' he is referring to here appear in his
book \cite[pp.\;119--120]{GG70} (rather than in Cauchy's article).
They are all formulated in terms of quantifier alternations, while the
variable~$x$ ranges through the real domain of the function.
Grattan-Guinness then paraphrases the revised version of Cauchy's
theorem from \cite[p.\;33]{Ca53}, asserts that the proof is based on
``the combined smallness of~$r_n(x_0)$ and~$r_n(x_0+\alpha)$''
\cite[p.\;123]{GG70} among other ingredients, and concludes that
``uniform convergence itself was tucked away in the word `always' with
no reference to the variable at all.'' (ibid.)  We summarize his
position as follows:
\begin{enumerate}
\item[GG1.]
Cauchy exploits infinitesimals in his attempted proof of the sum
theorem.
\item[GG2.]
More specifically, Cauchy evaluates the remainder term~$r_n$ at an
input~$x_0+\alpha$ where~$\alpha$ is an infinitesimal.
\item[GG3.]
Cauchy's proof is based on the smallness of~$r_n(x_0+\alpha)$.
\item[GG4.]  The procedures in \cite{Ca21} and \cite{Ca53} involve no
quantifier alternations.%
\footnote{Grattan-Guinness acknowledges this point GG4 implicitly in his
summary of Cauchy's proof in \cite[p.\;123]{GG70}.  Meanwhile,
an~$\varepsilon$ does occur in Cauchy on p.\;32; see
Section~\ref{42}.}
\item[GG5.]  
Cauchy's arguments relying on infinitesimals are difficult
to interpret in the context of the traditional characterisations of
uniform convergence which rely on quantifier alternations.
\item[GG6.]
Cauchy's term \emph{always} alludes to uniform convergence.
\item[GG7.]
Cauchy's term \emph{always} makes no reference to the variable~$x$.
\end{enumerate}

%
%

\subsection{Shortcomings}
\label{s42}

The analysis of Cauchy's procedures by Grattan-Guinness contains three
shortcomings.  First, claim~GG3 concerning the
expression~$r_n(x_0+\alpha)$ does not mention the crucial relation
\[
\alpha=\frac{1}{n}
\]
between~$n$ and the infinitesimal.  This relation is crucial to
Cauchy's analysis of series~(2) namely $\sin x+\frac{1}{2}\sin
2x+\frac{1}{3}\sin 3x+\ldots{}${} already discussed in
Section~\ref{23b}.  In fact in 1853 Cauchy did not use the
symbol~$\alpha$ at all but worked directly with the
infinitesimal~$\frac{1}{n}$.

As far as claim~GG6 is concerned, Grattan-Guinness makes no attempt to
explain the relation between uniform convergence and the term
\emph{always} exploited in the formulation of Cauchy's theorem.

Finally, his claim~GG7 is misleading.  The term \emph{always} is
indeed meant to refer to the variable~$x$ (as it did in Cauchy's
definition of continuity; see Section~\ref{toujours}).  Neither
Grabiner nor Grattan-Guinness properly address this point.  But in
order to explain precisely in what way the term refers to~$x$, one
needs to abandon the limitations of the particular interpretive
paradigm they are both laboring under, and avail oneself of a modern
framework where Cauchy's procedures are \emph{not} ``difficult to
interpret''; see Section~\ref{seven}.

\section{Giusti on the meaning of \emph{toujours}}
\label{five}

E. Giusti proposed two different interpretations of Cauchy's sum
theorem, which we explore in Sections~\ref{51b} and \ref{52}.

\subsection{Giusti's epsilontic track}
\label{51b}

In a 1984 study of Cauchy's sum theorem, E.\;Giusti wrote:
\begin{quote}
The comparison between the statement and the proof \emph{clarifies
immediately the meaning} of the word ``toujours'' in this context.
The sentence: ``the sum (3)%
\footnote{This is a reference to Cauchy's sum (3) namely
$u_n+u_{n+1}~\ldots u_{n'-1}$ discussed in Section~\ref{two}.}
always becomes infinitesimal for infinite~$n$ and~$n'> n$'' turns out
to mean nothing but that ``it is possible to assign to~$n$ a value
sufficiently large so that it follows that~$|s_{n'}-s_n|<\varepsilon$
for every~$n'>n$ and every~$x$''; \cite[p.\;38]{Gi84} (emphasis added)
\end{quote}
Giusti's comment indeed ``clarifies immediately the meaning'' of the
sum theorem but only in the sense of providing a mathematically
correct result loosely corresponding to Cauchy's theorem.  Giusti
continues:
\begin{quote}
in other words, it indicates the \emph{uniform} convergence of the
series.  As support for this, Cauchy shows how in the case of Abel's
series the condition is not satisfied, and the sum (3) does \emph{not}
become infinitely small for infinitely large values of~$n$
and~$n'$. (ibid.) (emphasis in the original; translation ours)
\end{quote}
To summarize, Giusti's formulation relies on thinly veiled quantifier
alternations, requires interpreting what Cauchy calls an ``infinite
number~$n$'' as a variable subordinate to a universal quantifier, and
sheds little light on Cauchy's term \emph{always}.  Cauchy's own
\emph{procedures} are therefore not fully clarified.

\subsection{Giusti's sequential track}
\label{52}

A few pages later, Giusti offers an alternative interpretation in the
following terms:
\begin{quote}
[In Cauchy's sum theorem] the magnitude which plays a crucial role is
the remainder term~$r_n(x)$, depending on two variables.  The
convergence of that series implies that~$r_n(x)$ always tends to~$0$
for infinite values of~$x$ [sic], which means, from the standpoint of
our interpretation, that for each sequence~$x_n$ the variable
$r_n(x_n)$ is infinitesimal.  Once again, this [procedure] is
equivalent to [proving] uniform convergence. \ldots{} [In Cauchy's
analysis of the series of Abel], as we have already explained, Cauchy
proves that if the sum~$u_n+u_{n+1}+\ldots+u_{n'-1}$ becomes always
infinitely small for some infinitely big values of~$n$ and~$n'$, then
the series~$u_1+u_2+u_3+\ldots$ converges and its sum is a continuous
function.  \cite[p.\;50]{Gi84} [translation ours]
\end{quote}
Thus, Giusti provides a sequential interpretation of Cauchy's
argument, still in an Archimedean context (see Section~\ref{64}, item
[C]).  However, as Laugwitz pointed out,
\begin{quote}
Giusti gives a correct translation of the example into the language of
sequences, \ldots{} [Giusti 1984, 50].  But he fails to translate the
general theorem and its proof.  Actually both~$x$ (or~$x + \alpha$)
and~$n$ will have to be replaced by sequences which becomes
troublesome as soon as they are not connected as in the example where
$x = 1/n$.  Moreover, the theorem shows the power of Cauchy's ``direct
consideration of infinitesimals.''  \cite[p.\;266]{La87}
\end{quote}
A paraphrase of Cauchy's proof along the lines of Giusti's sequential
reading may be harder to implement than the analysis of Abel's
example, and may depart significantly from Cauchy's own proof as
summarized in \cite[formula~(6.30), p.\;123]{GG70}, which we present
in Section~\ref{54}.  An additional difficulty (with Giusti's
[C]-track reading; see Section~\ref{64}) not mentioned by Laugwitz is
detailed in Section~\ref{order}.

\subsection{Order~$n'>n$}
\label{order}

If~$n$ and~$n'$ are (infinite) \emph{numbers} then one naturally
expects as Cauchy does that they are comparable; as Cauchy writes, we
have~$n'>n$.  However if they are \emph{sequences}~$n =(n_k)$ and~$n'
= (n'_k)$ tending to infinity, then there is no reason to assume that
they are comparable.  For example, if~$n_k=k$ whereas~$n'_k=k+(-1)^k$
then who is to say which sequence is bigger,~$n'$ or~$n$?  Meanwhile
in any ordered Bernoullian continuum (see Section~\ref{53})
comparability of numbers is automatic.  This is not to say that Cauchy
\emph{constructed} an ordered Bernoullian continuum, but rather that
his \emph{procedures} presume that the infinite numbers are ordered,
and this idea finds more faithful proxies in modern infinitesimal
theories than in Giusti's [C]-track reading.

\subsection{Cauchy's track}
\label{54}

To show that~$s(x_0+\alpha)-s(x_0)$ is infinitesimal, Cauchy writes
$r_n=s-s_n$, so that
\[
s(x_0+\alpha)-s(x_0)=
\left[s_n^{\phantom{I}}(x_0+\alpha)-s_n(x_0)\right] +
\left[r_n^{\phantom{I}}(x_0+\alpha)-r_n(x_0)\right]
\]
(see Section~\ref{s54}).  The terms~$r_n^{\phantom{I}}(x_0+\alpha)$
and~$r_n(x_0)$ are infinitesimal for \emph{each} infinite~$n$ by
Cauchy's \emph{always} hypothesis.  To account for the smallness of
the first summand, Grattan-Guinness appeals to the ``continuity
of~$s_n(x)$.''  However, this is immediate only for finite~$n$,
whereas Cauchy's proof of continuity requires~$\Delta
s_n=s_n^{\phantom{I}}(x_0+\alpha)-s_n(x_0)$ to be infinitesimal for
\emph{some} infinite~$n$.

Can one guarantee that if~$\Delta s_n$ is infinitesimal for each
finite~$n$, this property will persist for some infinite~$n$?  In the
context of hyperreal proxies for Cauchy's procedures, the answer is
affirmative (every internal sequence which is infinitesimal for all
finite~$n$ will remain infinitesimal for \emph{some} infinite~$n$).
What is perhaps even more remarkable than the validity of this type of
\emph{permanence principle} is the fact that a related principle can
indeed be found in \cite{Ca29}, as reported in \cite[p.\;274]{Ro66}.
The relevant technical result in a hyperreal framework is
\cite[Theorem~3.3.20, p.\;65]{Ro66}.  \cite[p.\;49]{MT97} follow
\cite{La90} in interpreting Euler in terms of \emph{hidden lemmas}
using a related result called the \emph{overspill theorem}; see also
the \emph{sequential theorem} in \cite[p.\;361]{MT01}.

\section{Conclusion}

Interpretation of texts written in the nineteenth century, and the
meaning we give to technical terms, procedures, theories, and the like
are closely related to what we already know as well as our
expectations and assumptions.  This paper provides evidence that a
change in the cultural-technical framework of a historian provides new
explanations, which are arguably more natural, and new insights into
Cauchy's work.

\section*{Acknowledgments} 
V.\;Kanovei was supported in part by the RFBR grant
number~17-01-00705.  M.\;Katz was partially funded by the Israel
Science Foundation grant number~1517/12.  We are grateful to Dave
L.\;Renfro for helpful suggestions.

\medskip

\textbf{Tiziana Bascelli} graduated in mathematics (1993) and
philosophy (2006), and obtained a PhD in theoretical philosophy (2010)
from University of Padua, Italy.  She co-authored \emph{Campano da
Novara's Equatorium Planetarum.  Transcription, Italian translation
and commentary} (Padua, 2007) and \emph{Gali\-leo's `Sidereus Nuncius'
or `A Sidereal Message'}, translated from the Latin by William
R. Shea, introduction and notes by William R. Shea and Tiziana
Bascelli.  Sagamore Beach, MA (USA): Science History Publications,
2009.  She is an independent researcher in integrated history and
philosophy of science, in the field of early modern mechanics and
mathematics.  Her research interests are in the development of
infinitesimal objects and procedures in Seventeenth-century
mathematics.

\medskip

\textbf{Piotr B\l aszczyk} is Professor at the Institute of
Mathematics, Pedagogical University (Cracow, Poland).  He obtained
degrees in mathematics (1986) and philosophy (1994) from Jagiellonian
University (Cracow, Poland), and a PhD in ontology (2002) from
Jagiellonian University.  He authored \emph{Philosophical Analysis of
Richard Dedekind's memoir \emph{Stetigkeit und irrationale Zahlen}}
(2008, Habilitationsschrift).  He co-authored \emph{Euclid,
\emph{Elements, Books V--VI}.  Translation and commentary}, 2013; and
\emph{Descartes, Geometry.  Translation and commentary} (Cracow,
2015).  His research interest is in the idea of continuum and
continuity from Euclid to modern times.

\medskip

\textbf{Alexandre Borovik} is Professor of Pure Mathematics at the
University of Manchester, United Kingdom.  His principal research
interests are to algebra, model theory, and combinatorics.  He
co-authored monographs \emph{Groups of Finite Morley Rank}, 1994,
\emph{Coxeter Matroids}, 2003, and \emph{Simple Groups of Finite
Morley Rank}, 2008.  His latest work, on probabilistic recognition of
so-called black box groups, focuses on the boundary between finite and
infinite in mathematics.

\medskip

\textbf{Vladimir Kanovei} graduated in 1973 from Moscow State
University, and obtained a Ph.D. in physics and mathematics from
Moscow State University in 1976. In 1986, he became Doctor of Science
in physics and mathematics at Moscow Steklov Mathematical Institute
(MIAN).  He is currently Principal Researcher at the Institute for
Information Transmission Problems (IITP) and Professor at Moscow State
University of Railway Engineering (MIIT), Moscow, Russia.  Among his
publications is the book \emph{Borel equivalence relations. Structure
and classification}. University Lecture Series 44. American
Mathematical Society, Providence, RI, 2008.

\medskip

\textbf{Karin U. Katz} (B.A. Bryn Mawr College, '82); Ph.D. Indiana
University, '91) teaches mathematics at Bar Ilan University, Ramat
Gan, Israel.  Among her publications is the joint article ``Proofs and
retributions, or: why Sarah can't \emph{take} limits'' published in
\emph{Foundations of Science}.

\medskip

\textbf{Mikhail G. Katz} (BA Harvard '80; PhD Columbia '84) is
Professor of Mathematics at Bar Ilan University, Ramat Gan, Israel.
He is interested in Riemannian geometry, infinitesimals, debunking
mathematical history written by the victors, as well as in true
infinitesimal differential geometry; see \emph{Journal of Logic and
Analysis} \textbf{7}:5 (2015), 1--44 at
\url{http://www.logicandanalysis.com/index.php/jla/article/view/237}

\medskip

\textbf{Semen S. Kutateladze} was born in 1945 in Leningrad (now
St.~Petersburg).  He is a senior principal officer of the Sobolev
Institute of Mathematics in Novosibirsk and professor at Novosibirsk
State University.  He authored more than 20 books and 200 papers in
functional analysis, convex geometry, optimization, and nonstandard
and Boolean valued analysis.  He is a member of the editorial boards
of \emph{Siberian Mathematical Journal}, \emph{Journal of Applied and
Industrial Mathematics}, \emph{Positivity}, \emph{Mathematical Notes},
etc.

\medskip

\textbf{Thomas McGaffey} teaches mathematics at San Jacinto College in
Houston, US.

\medskip

\textbf{David M. Schaps} is Professor Emeritus of Classical Studies at
Bar Ilan University, Israel.  His books include \emph{Economic Rights
of Women in Ancient Greece}, \emph{The Invention of Coinage and the
Monetization of Ancient Greece}, and \emph{Handbook for Classical
Research}; among his articles are ``The Woman Least Mentioned:
Etiquette and Women's Names", ``What was Free about a Free Athenian
Woman?'', and ``Zeus the Wife-Beater''.

\medskip

\textbf{David Sherry} is Professor of Philosophy at Northern Arizona
University, in the tall, cool pines of the Colorado Plateau.  He has
research interests in philosophy of mathematics, especially applied
mathematics and non-standard analysis.  Recent publications include
``Fields and the Intelligibility of Contact Action,'' \emph{Philosophy
90} (2015), 457--478.  ``Leibniz's Infinitesimals: Their Fictionality,
their Modern Implementations, and their Foes from Berkeley to Russell
and Beyond,'' with Mikhail Katz, \emph{Erkenntnis 78} (2013), 571-625.
``Infinitesimals, Imaginaries, Ideals, and Fictions,'' with Mikhail
Katz, \emph{Studia Leibnitiana 44} (2012), 166--192.  ``Thermoscopes,
Thermometers, and the Foundations of Measurement,'' \emph{Studies in
History and Philosophy of Science 24} (2011), 509--524.  ``Reason,
Habit, and Applied Mathematics,'' \emph{Hume Studies 35} (2009),
57--85.

\vfill\eject \end{document}